\documentclass[reqno,centertags]{amsart}[12pt]
\usepackage{amsmath,amsthm,amscd,amssymb}
\usepackage{latexsym}
\usepackage{upref}

\newcommand{\Pk}{{\Pi_\kappa}}
\newcommand{\dA}{{A}}
\newcommand{\bbA}{{\mathbb{A}}}
\newcommand{\bbR}{{\mathbb{R}}}

\newcommand{\bbC}{{\mathbb{C}}}

\newcommand{\DdA}{\calD(\dA)}

\newcommand{\mN}{\mathfrak N}
\newcommand{\mL}{\mathfrak L}
\newcommand{\mM}{\mathfrak M}

\newcommand{\calA}{{\mathcal A}}
\newcommand{\calB}{{\mathcal B}}

\newcommand{\calD}{{\mathcal D}}

\newcommand{\calH}{{\mathcal H}}

\newcommand{\calJ}{{\mathcal J}}

\newcommand{\calL}{{\mathcal L}}
\newcommand{\calM}{{\mathcal M}}
\newcommand{\calN}{{\mathfrak N}}

\newcommand{\calR}{{\mathcal R}}

\newcommand{\mh}{{\mathfrak H}}

\newcommand{\ol}{\overline}
\newcommand{\ti}{\tilde  }

\newcommand{\loc}{\text{\rm{loc}}}

\newcommand{\bi}{\bibitem}

\def\mul{{\rm mul\,}}
\newcommand{\range}{\mathfrak R}
\renewcommand{\Re}{\text{\rm Re}}
\renewcommand{\Im}{\text{\rm Im}}

\allowdisplaybreaks

\newtheorem{theorem}{Theorem}

\newtheorem{corollary}[theorem]{Corollary}

\theoremstyle{definition}
\newtheorem{definition}[theorem]{Definition}
\newtheorem{example}[theorem]{Example}

\theoremstyle{remark}
\newtheorem{remark}[theorem]{Remark}

\begin{document}
\title[On realization of the Krein-Langer class $N_\kappa$]{On realization
of the Krein-Langer class $N_\kappa$ of matrix-valued functions in Hilbert
spaces with indefinite metric}
\author[Yu.~Arlinskii, S.~Belyi, V.~Derkach, and E.~Tsekanovskii]{Yuri Arlinskii, Sergey Belyi, Vladimir Derkach,
and Eduard Tsekanovskii}
\address {Department of Mathematics \\
East Ukrainian National University \\
Kvartal Molodyozhny, 20-A \\
91034 Lugansk \\
Ukraine}
\email{yma@snu.edu.ua}
\address{Department of Mathematics\\ Troy State University\\
Troy, AL 36082, USA} \email{sbelyi@trojan.troyst.edu\\
\indent{\it URL:}http://spectrum.troyst.edu/$\sim$belyi}
\address{Department of Mathematics\\ Donetsk National University,
Universitetskaya str, 24, 83055 Donetsk, Ukraine}
 \email{derkach@univ.donetsk.ua}

%
\address{Department of Mathematics\\ Niagara University, NY
14109, USA} \email{tsekanov@niagara.edu\\
\indent{\it URL:}http://faculty.niagara.edu/tsekanov/}
\dedicatory{Dedicated to Henk de Snoo on the occasion of his 60th
birthday}

\date{\today}

\begin{abstract}
In this paper the realization problems for the Krein-Langer class
$N_\kappa$ of matrix-valued functions are being considered. We
found the criterion when a given matrix-valued function from the
class $N_\kappa$ can be realized as linear-fractional
transformation of the transfer function of canonical conservative
system of the M.~Livsic type (Brodskii-Livsic rigged operator
colligation) with the main operator acting on a rigged Pontryagin
space $\Pk$ with indefinite metric. We specify three subclasses of
the class $N_\kappa(R)$ of all realizable matrix-valued functions
that correspond to different properties  of a realizing system, in
particular, when the domains of the main operator of a system and
its conjugate coincide, when the domain of the hermitian part of a
main operator is dense in $\Pi\kappa$. Alternatively we show that
the class $N_\kappa(R)$ can be realized as  transfer
matrix-functions of some canonical impedance systems with
self-adjoint main operators in rigged spaces $\Pk$. The case of
scalar  functions of the class $N_\kappa(R)$ is considered in
details and some examples are presented.

 \end{abstract}

\maketitle
\section{Introduction}\label{s1}

Realizations and corresponding operator models of different
classes of holomorphic matrix-valued functions in the open right
half-plane, unit circle and upper half-plane play important role
in spectral analysis of  different classes of linear operators in
Hilbert spaces,
interpolation problems and system theory, and we refer in this
matter to \cite{ADRS}, \cite{ABDS}, \cite{ArDer}-\cite{Aro1},
\cite{BallSt}, \cite{BT2}-\cite{BT4}, \cite{Br}, \cite{DerHas},
\cite{DHS}, \cite{DTs}, \cite{GMKT}, \cite{KL1},
\cite{St1}-\cite{TSh1}. In this paper we continue the
investigation of various problems
 that arise in the study of linear stationary conservative dynamic
 systems (operator colligations). Relying on the results and
technique developed in \cite{BT2}, \cite{BT3} we keep dealing with
linear stationary conservative dynamic systems (l.s.c.d.s) $\theta$
of the form
$$
\begin{cases}
(\bbA&-zI)=K\calJ\varphi_-\\
\varphi_+&=\varphi_--2iK^\ast x
\end{cases} \qquad
{\left(\Im \,\bbA=K\calJ K^\ast \right) } \qquad
$$
or
\[
\theta =\left(\begin{array}{ccc}
\bbA &K &\calJ\\
\mh^+\subset \Pk\subset \mh_- &{ } &E
\end{array}\right).
\]
In the system $\theta$  above $\mathbb  A$ is a bounded linear
operator acting from $\mathfrak H^+$ into $\mathfrak H_-$, where
$\mathfrak H^+\subset\Pk\subset\mathfrak H_-$ is a rigged
Pontryagin space,  $K$ is a linear bounded operator from a Hilbert
space  $E$ into $\mathfrak H_-$, $\calJ=\calJ^\ast=\calJ^{-1}$ is
acting in $E$, $\varphi_\pm\in E$, $\varphi_-$ is an input vector,
$\varphi_+$ is an output vector, and $x\in\mathfrak H^+$ is a
vector  of the inner state of the system $\theta$.  The
operator-valued function
$$W_\theta (z)=I-2iK^\ast (\mathbb
A-zI)^{-1}K\calJ \qquad (\varphi_+= W_\theta (z)\varphi_-), $$ is
the transfer operator-valued function of the system $\theta$. It
was shown in \cite{BT3} that a Herglotz-Nevanlinna matrix-valued
function $V(z)$ acting on a Hilbert space $E$  can be represented
and realized in the form
 \[
V(z)=i[W_\Theta(z)+I]^{-1}[W_\Theta(z)-I]\calJ
 =K^*(\dA_{R}-zI)^{-1}K,
\]
where $W_\Theta(z)$ is a transfer function  of some canonical
system $\Theta$,
\[
 \Theta =\begin{pmatrix}
          \dA        & K      &  \calJ  \\
          \mh_+\subset\mh\subset\mh_-  &   &E
         \end{pmatrix},
\]
if certain conditions on integral representation of $V(z)$ are
met.  Alternatively, one can realize a Herglotz-Nevanlinna
matrix-valued $V(z)$ as a transfer mapping of an impedance system
$\Delta$ of the form
\begin{equation}\label{imp1}
    \left\{%
\begin{array}{ll}
    (\mathbb D-z I)x=K\varphi_-, \\
    \varphi_+=K^*x, \\
\end{array}%
\right.
\end{equation}
where $\mathbb D$ is a self-adjoint operators acting from $\mh_+$
into $\mh_-$ (see \cite{BT3}, \cite{BT4}). In this case associated
transfer function is given by
\[
 V_\Delta(z) =K^*(\mathbb D-z I)^{-1}K.
\]
In this paper we study similar realization problems but utilize a
new type of realizing systems whose main operator is acting on a
rigged Pontryagin space $\Pk$. The set of realizable functions
appears to be a subclass of the well known Krein-Langer's class
$N_\kappa$ also known as generalized Nevanlinna functions. In
Section \ref{s6} we specify three subclasses of the class
$N_\kappa(R)$ of all realizable matrix-valued functions that yield
different properties of operators in the realizing systems. It is
worth mentioning that in the case when $\kappa=0$ all the
subclasses coincide with the similar subclasses of realizable
Herglotz-Nevanlinna functions described in \cite{BT3}, \cite{BT4}.
 Section \ref{s7} uses a factorization formula from \cite{DLLSh}
 to provide applications of $N_\kappa(R)$ realizations to the scalar case
 when $E=\bbC$ while establishing  a connection with the class $N(R)$ of realizable
Herglotz-Nevanlinna functions.
The paper is concluded with several examples.

\section{Operators in Pontryagin spaces $\Pi_\kappa$}\label{s2}

We start with the basic construction following some results from
the theory of operators in $\Pi_\kappa$ spaces \cite{AG93},
\cite{IK}, \cite{KL}, \cite{KL1}. Let $\Pi_\kappa$ be a Pontryagin
space \cite{IK}, i.e., a Hilbert space $\calH$ where along with
the usual scalar product $(x,y)$ there is an indefinite scalar
product
\begin{equation}\label{e1-1}
    [x,y]=(Jx,y),
\end{equation}
where $J=P_+-P_-$ is a bounded linear operator such that $J=J^*$,
$J^2=I$, and $P_+$ and $P_-$ are complementary orthoprojections,
$P_++P_-=I$. Putting $\Pi_\pm =P_\pm\Pi_\kappa$ we have
\[
    \Pi_\kappa=\Pi_+\boxplus\Pi_-,\quad \dim \Pi_-=\kappa.
\]
Here and below the direct orthogonal sum with respect to an
indefinite scalar product \eqref{e1-1} is denoted by $\boxplus$
and called $\pi$-\textit{orthogonal sum}. Similarly, the
$\pi$-\textit{orthogonal complement} of a lineal $L$ will be
denoted by $L^{[\perp]}$. The positive definite $(x,y)$ and
indefinite  $[x,y]$ scalar products are related by
\begin{equation*}\label{e1-3}
    \begin{aligned}
        (x,y)&=[x_+,y_+]-[x_-,y_-],\\
        [x,y]&=(x_+,y_+)-(x_-,y_-),\\
    \end{aligned}
\end{equation*}
where $x=x_++x_-$, $y=y_++y_-$, $x_+,y_+\in\Pi_+$, and
$x_-,y_-\in\Pi_-$.

The set of vectors $f\in L$ that are $\pi$-orthogonal to $L$, i.e.
$f[\perp] L$ is called \cite{IK} the \textit{isotropic part} of
the linear manifold $L$. If the isotropic part of $L$ has non-zero
elements we say that the scalar product $[\cdot,\cdot]$ is
\textit{degenerate} \cite{KL} on $L$. $L_+$ (respectively, $L_-$,
$L_0$) will denote the set of all $x\in \Pi_\kappa$ for which
$[x,x]>0$ (respectively, $[x,x]<0$, $[x,x]=0$) and is called
\textit{positive} (\textit{negative, neutral}) part of $L$. Every
subspace $\calL\in\Pi_\kappa$ can be decomposed into a direct sum
of $\pi$-orthogonal subspaces
\begin{equation*}\label{e1-4}
    \calL=\calL_+\boxplus\calL_0\boxplus\calL_-,
\end{equation*}
where $\calL_+$, $\calL_0$, and $\calL_-$ are, respectively,
positive, neutral, and negative subspaces, some of which may
degenerate into null subspaces. For a subspace $\calL$ above we
write $sign\,\calL=(l_+,l_0,l_-)$ where $l_\pm=\dim \calL_\pm$ and
$l_0=\dim \calL_0$ \cite{KL}.

Recall \cite{ADRS} that a \textit{linear relation} in $\Pi_\kappa$
is a subspace ${\calA}$ in $\Pi_\kappa\times\Pi_\kappa$. The
\textit{domain} of a linear relation ${\calA}$ is
\[
\calD({\calA})=\left\{f\in\Pi_\kappa:\left<f,f'\right>\in
\calA\quad\mbox{for some}\quad f'\in\Pi_\kappa\right\},
\]
and the range of ${\calA}$ is
\[
\calR({\calA})=\left\{f'\in\Pi_\kappa:\left<f,f'\right>\in
{\calA}\quad\mbox{for some}\quad f\in\Pi_\kappa\right\}.
\]
The subspace
\[
\mul{\calA}=\left\{g\in\Pi_\kappa:\left<0,g\right>\in{\calA}\right\}
\]
is called the \textit{multivalued part} of a linear relation
${\calA}$. A linear relation ${\calA}$ is the graph of a linear
operator in $\Pi_\kappa$ if and only if $\mul{\calA}=\{0\}.$

Let us associate with a linear operator $A$ in $\Pi_\kappa$ the
linear relation ${\calA}:=Gr(A)$, the graph of the operator $A$.

 For a linear
relation ${\calA}$ in $\Pi_\kappa$ its $\pi$-\textit{adjoint}
${\calA}^+$ is defined by
\[
{\calA}^+=\left\{\left<h,h'\right>\in\Pi_\kappa\times\Pi_\kappa:[f',h]=[f,h']\quad\mbox{for
all}\quad \left<f,f'\right>\in{\calA}\right\}.
\]
A linear relation ${\calA}$ (operator $A$) is called
$\pi$-\textit{symmetric} if ${\calA}\subset{\calA}^+$  and
$\pi$-\textit{selfadjoint} if ${\calA}={\calA}^+$.

We recall \cite{IK} that a $\pi$-symmetric operator $\dA$ in $\Pk$
can not have more than $\kappa$ eigenvalues, counting
multiplicities, in the upper (lower) half-plane. If the operator $A$
is $\pi$-self-adjoint, then  these non-real eigenvalues are located
symmetrically with respect to the real axis. For an arbitrary
complex number $z$ and a $\pi$-symmetric operator $\dA$ in $\Pk$ we
set \cite{IK}
\[
    \calM_z=(\dA-z)\DdA,\qquad \calN_{\bar z}=\calM_z^{[\perp]}.
\]
 If $\lambda$ (Im$\lambda\ne 0$) is not an
eigenvalue of $\dA$, then $\calM_\lambda$ is a subspace of $\Pk$
and $\calN_\lambda$ is called \cite{IK} a \textit{deficiency
subspace} corresponding to $\lambda$ with $\dim \calN_\lambda$
maintaining a constant value as a \textit{deficiency index} of
$\dA$ in $\Pk$. Let $\Delta_{\dA}$ be the set of all non-real
$\lambda$ for which the scalar product $[\cdot,\cdot]$ is
degenerate on $\calN_\lambda$. According to \cite{KL} the set
$\Delta_{\dA}$ of a $\pi$-symmetric operator $\dA$ contains no
interior points, its complement
$(\bbC_+\cup\bbC_-)\setminus\Delta_{\dA}$ is an open set, and on
every component of this open set $sign\,\calN_\lambda$ is
constant.

It was shown in \cite{Der2} that every $\pi$-symmetric operator
$\dA$ in the space $\Pk$ admits $\pi$-self-adjoint extensions in
$\Pk$ if and only if its deficiency indices coincide.
An operator $\dA$ is called
\textit{prime} if it has no non-real eigenvalues and
\begin{equation}\label{Simple}
    \texttt{c.l.s.}\left\{\mN_z,\,z\ne\bar z\right\}=\Pk.
\end{equation}

In what follows we denote $\Re(T)=(T+T^+)/2$, $\Im(T)=(T-T^+)/2i$
for linear operators $T$  in $\Pk$ with $\calD (T) =\calD (T^+)$.
Similarly, for a linear operator $Q$ with $\calD (Q) =\calD (Q^*)$
in a Hilbert space we use the same notation to denote
$\Re(Q)=(Q+Q^*)/2$ and $\Im(Q)=(Q-Q^*)/2i$.

\section{Bi-extensions in Rigged Pontryagin Space}\label{s3}

 Let consider $\dA$ as an operator from the Hilbert space
 $\ol{\calD(\dA)}$ into the Hilbert space $\calH$. Then its
 adjoint $\dA^*$ is defined on a set $\calD(\dA^*)$ that is dense
 in $\calH$ and has the range in $\ol{\calD(\dA)}$. This allows us
 to introduce the Hilbert spaces $\mh_+=\calD(\dA^*)$ and
 $\mh^+=J\calD(\dA^*)$ with corresponding inner products
\[
\begin{split}
  &(f,g)_+ = (f,g)+(\dA^*f,\dA^*g),\quad f,\,g\in\mh_+, \\
  &(f,g)^+  = (f,g)+(\dA^*Jf,\dA^*Jg),\quad f,\,g\in\mh^+.
\end{split}
\]
Next we construct two 
rigged Hilbert spaces \cite{Ber}, \cite{Der}, \cite{ArDer},
\cite{BT3}
$$\mh_+\subset\calH\subset\mh_-\textrm{\quad and\quad
}\mh^+\subset\calH\subset\mh^-.$$ Let $\calR_1\in[\mh_+,\mh_-]$ and
$\calR_2\in[\mh^+,\mh^-]$ be the isometric Riesz-Berezanski\u{\i}
operators \cite{Ber} corresponding to the above triplets. We
introduce
$$J_+=J\Big|_{\mh_+},$$
the linear operator mapping $\mh_+$ isometrically onto $\mh^+$ and
let $J_+^\times\in[\mh^-,\mh_-]$ be its dual, which isometrically
maps $\mh^-$ onto $\mh_-$. It is easy to see that
$$\left(J_+^\times\right)^{-1}=\calR_2\,J_+\,\calR^{-1}_1.$$
Since $\mh_-$ is isomorphic to $\mh^-$ it can be considered as the
space of anti-linear functionals on $\mh^+$ defined by
$$\alpha(f)=(\alpha,Jf)=[\alpha,f],\quad
\alpha\in\mh_-,\,f\in\mh^+.$$ Thus, we can form a rigged $\Pk$
space
\[
    \mh^+\subset\Pk\subset\mh_-.
\]
Consequently, if $\bbA\in[\mh^+,\mh_-]$ then
$\bbA^\times\in[\mh^+,\mh_-]$ and $[\bbA f,g]=[f,\bbA^\times g]$
for all $f,g\in\mh^+$.
\begin{definition}\label{d1}
An operator $\bbA\in[\mh^+,\mh_-]$ is called \textit{bi-extension}
of a closed $\pi$-symmetric operator $\dA$ if $\bbA\supset\dA$ and
$\bbA^\times\supset\dA$. A bi-extension $\bbA$ is called
\textit{$\pi$-self-adjoint} if $\bbA=\bbA^\times$.
\end{definition}
Let $\bbA$ be a bi-extension of a $\pi$-symmetric operator $\dA$.
The operator $\hat A=\bbA$
with $\calD(\hat
A)=\left\{f\in\mh^+\mid \bbA f\in\Pk\right\}$ is called
\textit{quasi-kernel} \cite{ArDer} of the operator $\bbA$. If $\mathbb  A=\mathbb  A^\times$
and $\hat A$ is a quasi-kernel of $\mathbb  A$ such that $A\ne\hat
A$, $\hat A^+=\hat A$ then $\mathbb A$ is said to be a {\it
strong} $\pi$-self-adjoint bi-extension of $A$.

\textit{In what follows we will assume that the closed and
non-densely defined $\pi$-symmetric operator $\dA$ is
\textit{$J$-regular} \cite{ArDer}}, i.e., the operator $P\dA$ is
closed, where $P$ is the orthogonal projection onto
$\ol{J\calD(\dA)}$ in $\calH$. The analog of von Neumann's formula
for the operator $J\dA$ (see \cite{TSh1})
\[
    \mh^+=\calD(\dA)\oplus\mN'_i\oplus\mN'_{-i}\oplus\mN
\]
holds in the space $\mh^+$, where $\mN'_{\pm i}$ are the
\textit{semi-deficiency subspaces} of the operator $J\dA$
\cite{Krasn}, i.e.
$$
\mN'_{\pm i}=\ol{J\calD(\dA)}\ominus(PJ\dA\mp iI)\calD(\dA)),
$$
$\calR=J^\times_+\calR_2$, $\mN=\calR^{-1}\mL$, and
$\mL=\Pk\boxminus\ol{\calD(\dA)}$. The condition of $\dA$ being
$J$-regular is  equivalent to the subspace $\mL$ being closed in
$\mh_-$ and a sufficient condition of $J$-regularity is
$\dim\mL<\infty$ (see \cite{TSh1}, \cite{ArDer}). Let $P^+_{\DdA}$,
$P^+_{\mN'_i}$, $P^+_{\mN'_{-i}}$, $P^+_{\mN}$, and $P^+_{\mM}$ be
the orthogonal projections in $\mh^+$ onto $\DdA$, $\mN'_i$,
$\mN'_{-i}$, $\mN$, and $\mM=\mh^+\ominus\DdA$, respectively. Then
the set of all bi-extensions of a $J$-regular operator $\dA$ is
described by the formula \cite{TSh1}
\[
    \bbA=\dA P^+_{\DdA}+\left(\dA^++\calR S_\bbA\right)P^+_{\mM},
\]
where $S_\bbA\in[\mM,\mM]$. Moreover, $\bbA$ is $\pi$-self-adjoint
if and only if
\[
S_\bbA-S_\bbA^*=-iP^+_{\mN'_i}+iP^+_{\mN'_{-i}}.
\]
The Hilbert space version of the class $\Omega_{\dA}$ in the
definition below is found in \cite{Ar1}, \cite{Ar2}, \cite{TSh1}, \cite{BT2}.
\begin{definition}\label{d2}
We say that a closed densely defined linear operator $T$ acting in a
Pontryagin space $\Pi_\kappa$ belongs to the \textit{class}
$\Omega_A$ if:
\begin{enumerate}
\item[(1)] $T\supset \dA$, $T^+\supset \dA$ where $\dA$ is a
closed Hermitian operator;
\item[(2)] $T$ has a regular point in the lower half-plane;
\item[(3)] $PT$ and $PT^+$ are closed operators.
 \end{enumerate}
\end{definition}
 Note that a closed and non-densely defined $J$-regular $\pi$-symmetric
 operator $\dA$
 admits $\pi$-selfadjoint extensions of the class $\Omega_A$ if and only if its semi-deficiency indices coincide
 \cite{ArTsek}, \cite{TSh1}.

  An operator $\mathbb  A$ in  $[\mh^+,\mh_-]$ is called a
\textit{$(\ast)$-extension} \cite{BT2} of an operator $T$ of the
class $\Omega_A$  if both $\mathbb  A\supset T$ and
$\bbA^\times\supset T^+$. This $(\ast )$-extension is called {\it
correct} \cite{BT2} ({\it regular} \cite{Ar2}), if an operator
$\Re\,\mathbb  A=\frac{1}{2}(\mathbb  A+\mathbb  A^\times)$ is a
strong $\pi$-self-adjoint bi-extension of an operator $A$. It is
easy to show that if $\mathbb  A$ is a $(\ast )$-extension of $T$,
the $T$ and $T^+$ are quasi-kernels of $\mathbb  A$ and $\mathbb
A^+$, respectively.
\begin{definition}\label{d3} We say  the operator $T$ of the class
$\Omega_\dA$ belongs to the class $\Lambda_\dA$ if
\begin{enumerate}
\item[(1)] $T$ admits a correct $(\ast )$-extension; \item[(2)]
$Gr(A)=Gr(T)\cap Gr(T^+)$.
\end{enumerate}
\end{definition}
It can be shown (see \cite{ArDer}) that if $T\in\Lambda_\dA$ then the equation
$$(\bbA-\lambda I)x=g,$$
is solvable for all $\lambda\in\rho(T)$ and all $g\in \textrm{Im} \bbA
=\frac{1}{2}(\mathbb  A-\mathbb  A^\times)$.

\begin{remark}
\label{RR1} A survey of theory of bi-extensions of symmetric
operator in a Hilbert space and its application to characteristic
functions of operators of the class $\Lambda_A$ is presented in
\cite{TSh1}. Bi-extensions of $\pi$-symmetric operators in
Pontryagin spaces were studied in \cite{Der}.
\end{remark}
\section{Operator colligations in $\Pk$}\label{s4}
In this section we consider linear stationary conservative dynamic
systems (l. s. c. d. s.) $\theta $ of the form
\[
\begin{cases}
\aligned (\mathbb  A&-zI)=K\calJ\varphi_-\\
&\\
\varphi_+&=\varphi_--2iK^+ x
\endaligned
\end{cases}
\qquad {\left(\text{Im }\,\mathbb  A=K\calJ K^\ast \right). } \qquad
\]
In the system $\theta$  above $\mathbb  A$ is a bounded linear
operator acting from $\mathfrak H^+$ into $\mathfrak H_-$, where
$\mathfrak H^+\subset\Pk\subset\mathfrak H_-$ is a rigged
Pontryagin space,  $K$ is a linear bounded operator from a Hilbert
space  $E$ into $\mathfrak H_-$, $\calJ=\calJ^\ast=\calJ^{-1}$ is
acting in $E$, $\varphi_\pm\in E$, $\varphi_-$ is an input vector,
$\varphi_+$ is an output vector, and $x\in\mathfrak H^+$ is a
vector  of the inner state of the system $\theta$.

 For our purposes we
need the following more precise definition:
\begin{definition}\label{d3-1}
The array
\begin{equation}\label{e3-2}
\theta =\begin{pmatrix}
\mathbb  A &K &\calJ\\
\mathfrak H^+\subset\Pk\subset\mathfrak H_- &{ } &E
\end{pmatrix}
\end{equation}
is called a \textit{linear stationary conservative dynamic system}
(l.s.c.d.s.) or \textit{Brodski\u i-Liv\u sic rigged operator colligation}
if
\begin{enumerate}
\item[(1)] {$\mathbb  A$ is a correct ($\ast $)-extension of an
operator $T$ of the class $\Lambda_\dA$} for some $J$-regular
operator $\dA$ with finite and equal deficiency indices;
\item[(2)] {$\calJ=\calJ^\ast=\calJ^{-1}\in [E,E],\quad \dim E < \infty $};
\item[(3)] {$\bbA-\bbA^\times = 2iK\calJ K^+ $, where $K\in [E,\mh_-]
\quad \ker K=\{0\}\quad (K^+ \in [\mh^+,E])$}.
\end{enumerate}
\end{definition}

In this case, the operator $K$ is called a {\it channel operator}
and $\calJ$ is called a {\it direction operator} \cite{BT2}. We  associate with the system $\theta$ an operator-valued function
\[
W_\theta (z)=I-2iK^+ (\mathbb  A-zI)^{-1}K\calJ
\]
which is called {\it a transfer operator-valued function} of the
system $\theta$ or a \textit{characteristic operator-valued
function} of Brodski\u i-Liv\u sic rigged operator colligations
\cite{BT2}.

Following \cite{Br}, \cite{BT3} we call an l.s.c.d. system
$\theta$ \textsl{minimal} if the $\pi$-symmetric operator $A$ is
such that there are no nontrivial invariant subspaces on which $A$
induces  $\pi$-self-adjoint operators. Clearly, the l.s.c.d.
system $\theta$ is {minimal} if the operator $A$ satisfies the
condition~\eqref{Simple}.

 Let
$\theta$ be a l.s.c.d.s. of the form (\ref{e3-2}). We consider an
operator-valued function
\begin{equation}\label{e3-5}
V_\theta (z) = K^+ (\mathbb  A_R - zI)^{-1} K.
\end{equation}
The transfer operator-function $W_\theta (z)$ of the system $
\theta $ and an operator-function $V_\theta (z)$ of the form
(\ref{e3-5}) are connected by the  relation
\begin{equation}\label{e3-6}
V_\theta (z) = i [W_\theta (z) + I]^{-1} [W_\theta (z) - I] \calJ
\end{equation}

\section{Class $N_\kappa$. Realization Theorems.}\label{s5}

Let $E$ be a Hilbert space with an inner product $(\cdot,\cdot)$
and an operator-valued function $Q(z)$ belong to $[E,E]$.

\begin{definition} \mbox{\cite{KL1} } \label{d4}
An operator-valued function $V(z)\in[E,E]$ belongs to the class
$N_\kappa$ if it is meromorphic in $\bbC\setminus\bbR$ and such
that $V(\bar{z})={V(z)}^*$, $z \in Z_V$, and the kernel
    \begin{equation}
\label{kerN}
 {\mathsf{N}}_V(z,\zeta)
 =\frac{V(\zeta)-V(z)^*}{\zeta-\bar{z}},
 \quad z,\zeta \in Z_V,
 \quad
 \zeta\neq \bar{z},
\end{equation}
\[
 {\mathsf{N}}_V(z,\bar z)=V'(z),
 \quad z \in Z_V,
\]
has $\kappa$ negative squares, i.e. for all $z_j$ in the domain of
holomorphy $Z_V$ of the meromorphic  (in $\bbC\setminus\bbR$)
function $V(z)$ and $h_j\in E$ $( j=0,1,...,n,)$ the form
\begin{equation}\label{e6-43}
    \sum^n_{j,k=0}\Big(N_V(z_j,z_k)h_j,\,h_k\Big)\xi_j\,\bar\xi_k
\end{equation}
contains at most $\kappa$ negative squares and for one such a set
exactly $\kappa$ negative squares.
\end{definition}
Mention, that the kernel $N_V(z,\zeta)$ for a function $V\in
N_\kappa$ restricted to the upper half-plane has the same number
$\kappa$ of negative squares, see~\cite{KL}.

Class $N_\kappa$ was introduced  in \cite{KL1} and studied further
in \cite{KL2}, \cite{DL}. Different operator models corresponding
to $N_\kappa$-functions are constructed in \cite{KL1},
\cite{DerHas}, \cite{DHS}, \cite{DerHasSnoo}, \cite{GzZ}.

\begin{definition}\mbox{\cite{BT2} }\label{d4-1} An
operator-valued function $V(z)$ in a finite-dimensional Hilbert
space $E$ 
is called realizable if, in some domain ${\calD}\subset\bbC_-$,
$V(z)$ can be represented in the form
    \begin{equation}\label{Vtheta}
V_\theta(z) = i[W_{\theta }(z)+I]^{-1}[W_{\theta }(z)-I],
\end{equation}
where $W_\theta (z)$ is a transfer operator-function of some
l.s.c.d.s.  $\theta $ with the direction operator $\calJ=I$.
\end{definition}

\begin{definition}\label{d4-2} An operator-function $V(z)\in [E,E] \quad (\dim E<\infty )$
 belongs to the class $N_\kappa(R)$ if the following conditions are met:
\begin{enumerate}
\item $V\in N_\kappa$; \item  for all $f\in E$
    \begin{equation}\label{e6-57}
        \lim_{y\uparrow\infty}\frac{\big(V(iy)f,f\big)}{y}=0;
\end{equation}

\item 
For all $z\in Z_V$
\begin{equation}\label{e-23}
    \bigcap_{\zeta\in Z_V}\ker N_V(\zeta,z)=\{0\};
\end{equation}
\item for all $f\in\calB=\left\{f\in E\mid \lim_{y\uparrow\infty}
y\big(\Im\, V (iy)f,f\big)<\infty\right\}$
\begin{equation}\label{e6-59}
\lim_{y\uparrow\infty}  V (iy)f=0.
\end{equation}
   \end{enumerate}
\end{definition}
\begin{theorem} \label{t8}
Let $ \theta $ be a minimal l.s.c.d.s. of the form {\rm(\ref{e3-2})}
with $\calJ=I$ in $E$, $\dim E<\infty$. Then the operator-function
$V_\theta (z)$ of the form {\rm(\ref{Vtheta})} admits a holomorphic
continuation to a function $V(z)$ which belongs to the class
$N_\kappa(R)$.
\end{theorem}
\begin{proof}
For a l.s.c.d.s $\theta$ of the form \eqref{e3-2} consider
$$V (z) = K^+ (\mathbb  A_R - zI)^{-1} K,\quad z\in\rho(\hat A_R),$$
where 
$\hat A_R$ is the quasi-kernel of the operator 
$\bbA_R$.
As follows from~\eqref{e3-6} $V(z)$ is a
holomorphic continuation of the  function $V_\theta(z)$. We set
\[
    \Gamma_z=(\bbA_R-zI)^{-1}K, \quad z\in\rho(\hat A_R).
\]
It can be seen that the operator $\Gamma_z$ is invertible and the
following relation is valid
\[
  \Gamma_z=(\bbA_R-\zeta I)(\bbA_R-zI)^{-1}\Gamma_\zeta,\quad
  z,\quad z\in\rho(\hat A_R).
\]
Let $z_0\in\rho(\hat A_R)$ and let $\ti A$ be a $\pi$ - symmetric
operator defined as follows
\[
   \calD(\ti A)=\left\{f\in\calD(\hat A_R)\mid \left[(\hat
   A_R-\bar
   z_0I)f,\Gamma_{z_0}
    g\right]=0,\,\forall g\in E\right\},\; \ti A f=\hat A_Rf,\; f\in \calD(\hat A_R),
\]
Since
\begin{eqnarray}\label{e-27}
  \frac{V(z)-V^*(\zeta)}{z-\bar\zeta} &=& \frac{K^+(\bbA_R-zI)^{-1}K-K^+(\bbA_R-
  \bar\zeta I)^{-1}K}{z-\bar\zeta} \nonumber\\
   &=& K^+(\bbA_R-  \bar\zeta I)^{-1}(\bbA_R-zI)^{-1}K\nonumber\\
   &=&\Gamma_\zeta^+\Gamma_z,\quad z,\quad z\in\rho(\hat A_R),
\end{eqnarray}
the operator-function $V(z)$ is a Krein-Langer $Q$-function
\cite{KL}, \cite{KL1} for the $\pi$ - symmetric operator $\ti A$ and
its $\pi$-self-adjoint extension $\hat A_R$. 
Let $f\in\calD(\ti A)$. Then for any $g\in\mh^+$
$$[\Im\, \bbA f,g]= [f,\Im\, \bbA g]=[f, (\bbA_R-\bar z_0I)\Gamma_{z_0} K^+g]=
[(\hat A_R-\bar z_0I)f,\Gamma_{z_0} K^+g]=0.$$ Thus, $\Im\, \bbA
f=0$ for $f\in\calD(\ti A)$ and $\mathbb  A\supset T\supset\ti A$,
$\bbA^\times\supset T^+\supset\ti A$. But this is possible only
when $\ti A=\dA$ since $\dA$ is the maximal $\pi$-symmetric part
of the operator $T$. Consequently, ${\Gamma_z E}=\mN_z$ and
$V_\theta(z)$ is a Krein-Langer $Q$-function for a prime $\pi$ -
symmetric operator $A$. It was shown in \cite{KL1} that a
$Q$-function of a prime $\pi$ - symmetric operator (in the upper
half-plane) belongs to the class $N_\kappa$ and condition
\eqref{e6-57} holds.

In order to prove \eqref{e-23} we note that since the operator
$\dA$ is prime then for all $f\in\bigcap_{\zeta\in Z_V}\ker
N_V(\zeta,z)$ and all $g\in E$ we have
$$\left[\Gamma_z f,\Gamma_\zeta g\right]=\left(\frac{V_\theta(z)-V^*_\theta(\zeta)}
{z-\bar\zeta}f,g\right)=0,$$ and hence $f=0$.

It was shown in \cite{KL1} that
\begin{equation}\label{B}
\mN_z\cap\calD(\hat A_R)=\Gamma_z \calB,\quad  z\in\rho(\hat A_R).
\end{equation}
It is easy to see that
\begin{equation}\label{L}
(\bbA_R-zI)(\mN_z\cap\calD(\hat A_R))=\mL,
\end{equation}
 where $\mL=\calD(\dA)^{[\perp]}$. Clearly, $\dim\calB<\infty$, since  $\calB\subset E$.

In order to prove \eqref{e6-59} we need to use the spectral
decomposition of  the $\pi$-self-adjoint operator $\bbA_R$. It was
shown in \cite{ArDer} that for all $g\in\calD(\bbA_R)$
\begin{equation}\label{e-28}
\bbA_R
g=\lim_{z=iy\uparrow\infty}(-z^2)\left[(\bbA_R-zI)^{-1}+\frac{1}{z}\right]g.
\end{equation}
It follows from \eqref{e-27} that
\begin{equation}\label{e-30}
    (-z^2)\Gamma^+_{z_0}\left[(\bbA_R-zI)^{-1}+\frac{1}{z}\right]\Gamma_{z_0}=
    \frac{-z^2[V_\theta(z)-V_\theta^*(z_0)]}{(z-z_0)(z-\bar
    z_0)}+\frac{zz_0}{z-z_0}\Gamma^+_{z_0}\Gamma_{z_0}
\end{equation}
Taking the limit in \eqref{e-30} and using \eqref{e-28}, we obtain
for all $ f\in\calB$
\begin{align}
\Gamma^+_{z_0}\bbA_R
\Gamma_{z_0}f&=-\lim_{z=iy\uparrow\infty}V_\theta(z)f+V_\theta^*(z_0)f+z_0\Gamma^+_{z_0}\Gamma_{z_0}f
\label{e6-31}\\
&= -\lim_{z=iy\uparrow\infty}V_\theta(z)f+\Re\,
V_\theta(z_0)f+\Re\, z_0\Gamma^+_{z_0}\Gamma_{z_0}f.\nonumber
\end{align}
Taking into account that $\Re\, V_\theta(z_0)
=\Gamma^+_{z_0}\bbA_R\Gamma_{z_0}f-\Re\,z_0\Gamma^+_{z_0}\Gamma_{z_0}$
we get \eqref{e6-59}.
 \end{proof}
\begin{theorem} \label{t9}
Let an operator-valued function $V(z)$  belong to the class
$N_\kappa(R)$.
Then $V(z)$ admits a minimal  realization by a system $\theta$ of
the form \eqref{e3-2} with $\calJ=I$.
\end{theorem}
\begin{proof} 
It was shown in \cite{KL1} that if an
operator-function $V(z)$ in $[E,E]$ satisfies the conditions 1)
and 2) of Definition~\ref{d4-2} then it is possible to construct a
Pontryagin space $\Pk(V)$ and a symmetric operator $\dA$ with a
$\pi$-self-adjoint extension $A_V$ in $\Pk(V)$.
The strictness condition 3) 
guarantees that the function $V(z)$ is a Krein-Langer $Q$-function
of a pair $\dA$ and $A_V$. We will use a reproducing kernel space
model for the operators $A$, $A_V$ elaborated in \cite{ABDS},
\cite{DM} and \cite{DHS}.

{\it Step 1}. Let us consider the reproducing kernel Pontryagin
space $\Pk(V)$ corresponding to the kernel $N(z, \zeta):=N_V(z,
\zeta)$. The latter means see \cite{ADRS} that $\Pk(V)$ consists
of functions holomorphic on $Z_V$ and for each $z\in Z_V$ and
$h\in E$ the followings hold:
\begin{enumerate}
\item[(a)] $N(z,\zeta)h$ belongs to $\Pk(V)$ as a function on
$\zeta$; \item[(b)] $[f(\cdot),N(z,\zeta)h]=(f(z),h)_E$ for every
$f(\cdot)$ in $\Pk(V)$.
\end{enumerate}
In particular, it follows from (b) that the evaluation operator
$f\mapsto (f(z),h)_E$ is continuous in $\Pk(V)$. The
multiplication operator
\begin{equation}\label{Amult}
A:f(\zeta)\mapsto \zeta f(\zeta),
\end{equation}
with the domain $\mbox{dom }A=\{f\in\Pk(V): \zeta f(\zeta)\in
\Pk(V)\}$ is a prime closed symmetric operator in $\Pk(V)$. Let
$\calA:=Gr(A)$ As was shown in~\cite{DHS} the adjoint linear
relation ${\calA}^+$ takes the form
\[
  \calA^+=\{\, \{f,\widetilde f\}\in \Pk(V)^2:\,
  \widetilde f(\zeta)-\zeta f(\zeta)=h_1-V(\zeta)h_0;
   \ h_0, h_1 \in E \,\}
\]
and the linear relation
\begin{equation}\label{AV}
\calA_V=\{\, \{f,\widetilde f\}\in \Pk(V)^2:\,
  \widetilde f(\zeta)-\zeta f(\zeta)=h_1\in E\,\}
\end{equation}
is a self-adjoint extension of $A$ with $\rho(A_V)=Z_V$.

The  deficiency subspace $\mN_z$ of $A$ $(z\in Z_V)$ consists of
vector-functions $\gamma(z)h:=N(\bar z,\cdot)h$, $h\in E$. The
mapping $\gamma(z):E\to\mN_z$ is injective since the assumption
$\gamma(z)h=0$ and the equality
\begin{equation}\label{Qfun}
    \left[N(\bar z,\cdot)h,N(\bar w,\cdot)g  \right]=\left(N(\bar z,\bar
    w)h,g\right)_E=\left(\frac{V(z)-V(w)^*}{z-\bar w}h,g\right)_E
\end{equation}
imply that $h\in\bigcap\limits_{w\in Z_V}\ker(V(z)-V(w)^*)$. It
follows from the hypothesis 3) that $h=0$. The linear span of
deficiency subspaces $\mN_z$  $(z\in Z_V)$ is dense in
$\Pi_\kappa(V)$ and, hence, the operator $V$ is a prime symmetric
operator in $\Pi_\kappa(V)$.

Let us show that $\gamma(z)$ satisfies the identity
\begin{equation}\label{Gfield}
    \gamma(z)=(\calA_V-z_0)(\calA_V-z)^{-1}\gamma(z_0), \quad z,z_0\in Z_V.
\end{equation}
This is straightforward from the identities
\begin{equation}\label{1}
    zN(\bar z,\zeta)h-\zeta N(\bar z,\zeta)h=V(z)h-V(\zeta)h,
\end{equation}
\begin{equation}\label{2}
    z_0N(\bar z_0,\zeta)h-\zeta N(\bar z_0,\zeta)h=V(z_0)h-V(\zeta)h.
\end{equation}
Subtracting \eqref{2} from \eqref{1} one obtains due to~\eqref{AV}
\begin{equation}\label{AVN}
   \left\{N(\bar z,\zeta)h-N(\bar z_0,\zeta),zN(\bar z,\zeta)h-z_0
    N(\bar z_0,\zeta)h\right\}\in \calA_V.
\end{equation}
The latter equality is equivalent to~\eqref{Gfield}.

The identities~\eqref{Gfield} and~\eqref{Qfun} show that $V(z)$ is
the Kre\u{\i}n-Langer $Q$-function of the pair $\calA$, $\calA_V$.
The assumption \eqref{e6-57} implies that $\mbox{mul
}\calA_V=\{0\}$, see~\cite{KL}, and, therefore, $\calA_V$ is the
graph of an operator $A_V$.

 {\it
Step 2}. The operator $A$ need not be densely defined. Let us
calculate vector-functions from the subspace
$\calL=\calD(A)^{[\perp]}$ explicitly. As follows from~\eqref{B}
and\eqref{L}
\begin{equation}\label{3}
    \gamma(z)\calB=\mN_z\cap \calD(A_V)=(A_V-z)^{-1}\calL.
\end{equation}
This implies that for every $h\in\calB$ there exists a strong
limit of the vector-function
\begin{equation}\label{4}
    h_\infty(\cdot):=\lim_{y\uparrow\infty}(iy)\gamma(-iy)h=\lim_{y\uparrow\infty}(iy)N(iy,\cdot)h
\end{equation}
as $y\to\infty$. Since the evaluation operator is continuous in
$\Pk(V)$ one obtains from the hypothesis 4) for every $z\in Z_V$
\begin{equation}\label{5}
\begin{split}
h_\infty(z)&=\lim_{y\uparrow\infty}(iy)N(iy,z)h\\
&=\lim_{y\uparrow\infty}(iy)\frac{V(z)-V(-iy)}{z+iy}h=V(z)h.
\end{split}
\end{equation}
Therefore, the subspace $\calL$ takes the form
\[
\calL=\left\{V(\cdot)h:\,h\in\calB\right\}.
\]
It follows from~\eqref{1} and \eqref{AV}  that for every $h\in
\calB$
\begin{equation}\label{6}
    (A_V-z)^{-1}V(\cdot)h=N(\bar z,\cdot)h\quad (z\in Z_V).
\end{equation}
One can derive the same equality from~\eqref{3}, ~\eqref{4},
~\eqref{5} and ~\eqref{Gfield}.

{\it Step 3}.
 Let us show that for every $g\in E$ the function
$V(\cdot)g$ generates a functional on
$\calD(A^+)=\calD(A_V)+\mN_i$ by the formulas
\begin{equation}\label{7}
    \left[(A_V-z_0)^{-1}f(\cdot),V(\cdot)g\right]=\left(f(z_0),g\right)_E, \quad
    f\in\Pk(V),
\end{equation}
\begin{equation}\label{8}
    \left[N(\bar z_0,\cdot)h,V(\cdot)g\right]=(V(z_0)h,g)_E, \quad
    h\in E,
\end{equation}
which is continuous in the norm of $\mh^+=\calD(A^+)$.

Mention first that the formulas~\eqref{7} and~\eqref{8} are
consistent since for $f(\cdot)=V(\cdot)h$ $(h\in\calB)$ one
obtains from~\eqref{6} and \eqref{7} the formula
\[
\left[(A_V-z_0)^{-1}V(\cdot)h,V(\cdot)g\right]=(V(z_0)h,g)_E,
\]
which agrees with~\eqref{8}.

Next, it follows from \eqref{7}, \eqref{8}, \eqref{Gfield} and the
identity
\[
(A_V-z)^{-1}=(A_V-z_0)^{-1}(A_V-z_0)(A_V-z)^{-1}
\]
that
\begin{equation}\label{9}
    \left[N(\bar z,\cdot)h,V(\cdot)g\right]
    =\left(\left[V(z_0)v+(z-z_0)\frac{V(z)-V(z_0)}{z-z_0}h\right],g\right)_E
    =(V(z)h,g)_E.
\end{equation}

Let now a sequence $\varphi_n\in\calD(A_V)$ converges to
$\varphi\in\calD(A_V)$ in $\mh^+$-norm. Then
\[
f_n(\cdot)=(A_V-z_0)\varphi_n(\cdot)\to
f(\cdot)=(A_V-z_0)\varphi(\cdot)\mbox{ strongly in }\Pk(V)
\]
and by continuity of the evaluation operator one has
\[
f_n(z)\to f(z)\quad\forall z\in Z_V.
\]
Then it follows from~\eqref{7} that
\[
[\varphi_n(\cdot),V(\cdot)g]\to[\varphi(\cdot),V(\cdot)g]\quad
(n\to\infty)
\]
and, therefore, the functional generated by $V(\cdot)g$
via~\eqref{7} and \eqref{8} is continious, since
$\dim\mh^+(\mbox{mod }\calD(A_V))<\infty$.

{\it Step 4}. Using the operator $\dA$ defined in \eqref{Amult} we
 construct a rigged Pontryagin space $\mh^+\subset\Pk(V)\subset\mh_-$ the way it was
 described in the section \ref{s3}.
The functional $V(\cdot)h$, $h\in E$ considered above can be viewed
as an element from $\mh_-$.
 Let us define a linear operator $K:\,E\to\mh_-$ by the
equality
\begin{equation}\label{10}
    Kh:=V(\cdot)h,\quad h\in E.
\end{equation}
Clearly the operator $K$ is invertible, otherwise $ V(z)\equiv 0$
and definition \ref{d4-2} is violated.   Let us extend the
operator $A_V$ to the linear operator
$\bbA_R:\,\mh^+=\calD(A_V)+\mN_{z_0}\to\mh_-$ by the equality
\begin{equation}\label{11}
    \bbA_RN(\bar z_0,\cdot)g=z_0N(\bar z_0,\cdot)g+V(\cdot)g,\quad g\in E.
\end{equation}
This definition agrees with~\eqref{6} for $g\in\calB\subset E$.
The operator $\bbA_R$ is a $\pi$-self-adjoint bi-extension of
$\dA$ with the quasi-kernel $A_V$. It follows from~\eqref{AVN}
and~\eqref{11} that
\begin{equation}\label{12}
    (\bbA_R-z) N(\bar z,\cdot)g=V(\cdot)g,\quad g\in E.
\end{equation}
Making use of~\eqref{10} and \eqref{12} one obtains
\[
K^+(\bbA_R-z)^{-1}Kg=K^+N(\bar z,\cdot)g,\quad z\in Z_V.
\]
An application of \eqref{9} and \eqref{10} implies
\[
(K^+N(\bar z,\cdot)g,h)=[N(\bar z,\cdot)g,V(\cdot)h]=(V(z)g,h)_E,
\]
and, therefore,
\[
K^+(\bbA_R-z)^{-1}K=V(z),\quad z\in Z_V.
\]

{\it Step 5}. We set $\bbA_I=KK^+$ and define
\[
    \bbA=\bbA_R+i\bbA_I.
\]
It is obvious that the l.s.c.d.s.
\begin{equation}\label{sys_theta}
    \theta=\begin{pmatrix}
\mathbb  A &K &I\\
\mathfrak H^+\subset\Pk(V)\subset\mathfrak H_- &{ } &E
\end{pmatrix}
\end{equation}
satisfies conditions (2) and (3) of Definition \ref{d3-1}. What
remains to show then is that the quasi-kernel $T$ of the operator
$\bbA$ belongs to the class $\Lambda_A$. To do this it suffices to
show that $\rho(T)\cap\bbC_-$ is nonempty and check that $\dA$ is
the maximal symmetric part of $T$ and $T^+$. Since $V\in N_\kappa$
it follows from \cite[Theorem~2.2]{KL2} that $V(z)-iI$ has at most
$\kappa$ zeros in $\bbC_-$. Let us assume without loss of generality
that the operator  $(V(z_0)-iI)$ is invertible. Consequently, the
operator
\[
H=I+i K^+(\bbA_R-z_0I)^{-1}K=I+i V(z_0)\in[E,E],
\]
is invertible as well. It follows from~\eqref{5}, ~\eqref{11},
~\eqref{12} that
\[
\bbA N(\bar z_0,\cdot)g=z_0N(\bar z_0,\cdot)g+V(\cdot)Hg,\quad g\in
E.
\]
Since $H$ is invertible this implies $\range(\bbA-z_0I)\supset
\range(K)$.

Let $f\in\mN_z$. Then $(\bbA-z_0I)f=g+(z-z_0)f$, where $g=(\bbA_R-z
I)f+i\bbA_If.$ Since $\range(\bbA_I)\subset{\range(K)}$ and
$(\bbA_R-z I)\mN_z={\range(K)}$, we have $g\in{\range(K)}$ and
therefore there is an $x$ such that $(\bbA-z_0I)x=g$. Thus
$$(\bbA-z_0I)(f-x)=(z-z_0)f,$$
i.e., $\range(\bbA-z_0I)\supset\mN_z$. Since the operator $\dA$ is
prime, $\texttt{c.l.s.}\left\{\mN_z,\,z\ne\bar z\right\}=\Pk(V)$,
we have $\ol{\range(T-z_0I)}=\Pk(V)$. On the other hand, the
relations
$$\range(T-z_0I)\supset\mM_{z_0},\quad
\dim\range(T-z_0I)(\mod(\mM_{z_0}))<\infty,$$ imply that
$\range(T-z_0I)$ is closed and hence coincides with $\Pk(V)$.
Similarly, one proves that $\range(T^+-z_0I)=\Pk(V)$,  $T^+$ is a
quasi-kernel of $\bbA^\times$, and concludes that $\bbA$ is a
correct $(*)$-extension of $T$.

To prove that $\dA$ is the maximal symmetric part of $T$ and
$T^+$, we assume the contrary. Then there exists such a
$\pi$-symmetric operator $A_0$ that
$$T\supset A_0\supset\dA,\quad T^+\supset A_0\supset\dA.$$
Consequently, for all $f\in\calD(A_0)$, $\bbA f=\bbA^\times f=A_0
f$,
$$\bbA_I f=\frac{(\bbA-\bbA^\times)}{2i}f=0,\quad
\bbA_R f=(\bbA-\bbA_I)f=A_0 f\in\Pk(V).$$ Thus, $f\in\calD(A_V)$
and $A_Vf=A_0 f$. For any $g\in\mh^+$ we have
$$[(A_V-\bar z_0I)f,\gamma(z_0) K^+g]=[f,\bbA_I g]=[\bbA_If,g]=0.$$
But since ${\range(\gamma(z_0) K^+)}=\mN_{z_0}$, then
$f\in\calD(\dA)$ and $\calD(A_0)=\calD(\dA)$. Thus $\dA$ is the
maximal symmetric part of $T$ and $T^+$.

 \end{proof}
\begin{remark}\label{r10}
We should mention that the realization results obtained in
Theorem~\ref{t9} can be interpreted as realization with impedance
systems $\Delta$ of the form \eqref{imp1} with $\mathbb D=\bbA_R$
and
$$V_\Delta(z)=K^+(\bbA_R-zI)^{-1}K. $$
\end{remark}
\begin{remark}\label{r11} In the recent paper \cite{RoSa04} the authors
derive an alternative integral representation for matrix-functions
of the class $N_\kappa$. It can be easily shown that all functions
of the class $N_\kappa(R)$ fall into the special class described
in Theorem 4.1 of \cite{RoSa04} and permit a reduced Krein-Langer
integral representation developed in that theorem.
\end{remark}

\section{Subclasses of the class $N_\kappa(R)$}\label{s6}

In this section we follow \cite{BT2} and introduce three distinct subclasses of the
class   of realizable operator-valued functions $N_\kappa(R)$.

 \begin{definition}\label{d10} An operator-function
$V(z)$ of the class $N_\kappa(R)$ belongs to the subclass
$N^0_\kappa(R)$ if in the definition \ref{d4-2} the subspace $\calB$
is trivial, i.e.,
    \begin{equation}\label{e-53}
        \lim_{y\uparrow\infty} y\big(\Im\, V (iy)f,f\big)=\infty,\quad \forall f\in E,\, f\ne0.
    \end{equation}
\end{definition}

\begin{definition}\label{d11} An operator-function
$V(z)$ of the class $N_\kappa(R)$ belongs to the subclass
$N^1_\kappa(R)$ if in the definition \ref{d4-2} the subspace $\calB=E$, i.e.,
\begin{equation}\label{e-54}
        \lim_{y\uparrow\infty} y\big(\Im\, V (iy)f,f\big)<\infty,\quad \forall f\in E,\, f\ne0.
    \end{equation}
\end{definition}

\begin{definition}\label{d12} An operator-function
$V(z)$ of the class $N_\kappa(R)$ belongs to the subclass
$N^{01}_\kappa(R)$ if in the definition \ref{d4-2} the subspace $\calB$ is neither trivial
nor equals $E$, i.e.
\[
    \{0\}\varsubsetneq\calB\varsubsetneq E.
\]
\end{definition}
One may notice that $N(R)$ is a union of three distinct subclasses
$N_0(R)$, $N_1(R)$ and $N_{01}(R)$. Now we prove the direct and
inverse   realization theorems in each of the subclasses.

\begin{theorem} \label{t10}
Let $ \theta $ be a l.s.c.d.s. of the form {\rm(\ref{e3-2})} such
that $\calJ=I$ and $\dA$ is an operator with dense domain. Then
operator-function $V_\theta (z)$ of the form \eqref{e3-5},
\eqref{e3-6} has a holomorphic continuation
 $V(z)$ which belongs to the class $N^0_\kappa(R)$.

 Conversely, let an operator-valued function $V(z)$  belong to the class
$N^0_\kappa(R)$.  Then $V(z)$ admits a minimal realization by a
system $\,\theta$ of the form \eqref{e3-2} with a densely defined
operator $\dA$.
\end{theorem}
\begin{proof}
Since the operator $A$ is densely defined $\calD(\widetilde A)\cap
\mN_z=\{0\}$ for every self-adjoint extension  $\widetilde A$ of
$A$. In particular, one obtains $\calD(\widehat A_R)\cap
\mN_z=\{0\}$. Due to~\eqref{B} this implies $\calB=\{0\}$.

Conversely, if $\calB=\{0\}$ then it follows from~\eqref{B} and
\eqref{L}  that $\calL=\{0\}$, and, hence, the operator $A$ is
densely defined.
\end{proof}
In order to proceed with the similar results in the class
$N^1_\kappa(R)$ we need to recall the definition of $O$-operator
\cite{BT3} and give its analogue for the spaces with indefinite
metric. A  $J$-regular $\pi$-symmetric operator $\dA$ is called an
$O$-operator if its semi-deficiency indices are equal to zero. As it
was shown in \cite{ArDer} for an operator $T\in\Lambda_\dA$ there
exist linear operators
\begin{eqnarray*}
    &M_T:\mN'_i\boxplus\mN\rightarrow\mN'_{-i}\boxplus\mN,\\
    &M_{T^+}:\mN'_{-i}\boxplus\mN\rightarrow\mN'_{i}\boxplus\mN,
\end{eqnarray*}
such that $\calM_T=Gr(M_T)$ and $\calM_{T^+}=Gr(M_{T^+})$ have
trivial intersections with the manifold $\{\{x,x\}:\,x\in\mN\}$ and
\begin{eqnarray}
  \calD(T) &=& \calD(\dA)\boxplus\{x-y\mid \langle x,y\rangle\in \calM_T\}, \label{e-57}\\
    \calD(T^+) &=& \calD(\dA)\boxplus\{x-y\mid \langle x,y\rangle\in
  \calM_{T^+}\}.\label{e-58}
\end{eqnarray}
 If $\dA$ is an $O$-operator then $M_T$ and $M_{T^+}$ are operators in
 $\mN$ ($\dim\mN<\infty$)
with $1\in\rho(M_T)\cap\rho(M_{T^+})$ and the relations
\eqref{e-57}-\eqref{e-58} imply that  $\calD(T)=\calD(T^+)=\mN$,
$\Im (T)$ is a bounded selfadjoint operator in $\Pi_\kappa$, and
$\Re(T)$ is a $\pi$-selfadjoint extension of $A$ in $\Pi_\kappa.$

\begin{theorem} \label{t15}
Let $ \theta $ be a l.s.c.d.s. of the form {\rm(\ref{e3-2})} such
that $\calJ=I$, $\dA$ is an $O$-operator, and $\calD(T)=\calD(T^+)$.
Then operator-function $V_\theta (z)$ of the form \eqref{e3-5},
\eqref{e3-6} has a holomorphic continuation
 $V(z)$ which belongs to the class $N^1_\kappa(R)$.

 Conversely, let an operator-valued function $V(z)$  belongs to the class
$N^1_\kappa(R)$.  Then $V(z)$ admits a minimal realization by a
system $\theta$ of the form \eqref{e3-2} with a non-densely defined
$O$-operator $A$.
\end{theorem}
\begin{proof} Once again using Theorem \ref{t8} we have that $V(z)\in N_\kappa(R)$ and
need to show  \eqref{e-54}. Since $\dA$ is an $O$-operator,
$\mN'_{z}=0$, and hence $\dim\mN_{z}=\dim\mN$. Using \eqref{B} we
have
$$
\dim(\Gamma_z\calB)=\dim(\calD(\widehat A_R)\cap \mN_z)=\dim\mN.
$$
On the other hand, $\mN_z=\Gamma_z E$ and thus $\dim\Gamma_z
E=\dim\mN$. Consequently, since $\calB$ is a subspace of $E$ and
$\dim E=\dim\calB$ we have $\calB=E$.

Conversely, if $\calB=E$, then using \eqref{B} we get
$\Gamma_z\calB=\Gamma_z E=\calD(\widehat A_R)\cap \mN_z=\mN_z$.
Applying \eqref{L} yields
$$
\dim\left((\bbA_R-zI)\mN_z\right)=\dim\mL=\dim\mN,\quad
z\in\rho(\hat A_R).
$$
Since $(\bbA_R-zI)$ is an invertible operator for $z\in\rho(\hat
A_R)$ we conclude that $\dim\mN_z=\dim\mN$. It can be shown (see
\cite{ArDer}, \cite{TSh1}) that the operator $P^+_\mM$ described
in the section \ref{s3} is a bijective mapping from $\mN_{\pm i}$
onto $\mN'_{\pm i}\boxplus\mN$. Considering the above we have then
$\dim\mN'_{\pm i}\boxplus\mN=\dim\mN$. This proves that $\mN'_{\pm
i}=0$ and thus $\dA$ is an $O$-operator.
\end{proof}

\begin{remark}
If $V\in N^0_\kappa(R)$ then the operator $A$ in the
realization~\eqref{e3-2} is densely defined and this implies that
$\calD(T)\ne\calD(T^+)$ for the operator $T$ mutually disjoint with
$T^+$. When the operator $A$  is nondensely defined even mutually
disjoint operators $T$ and $T^+$ may have the same domain. In fact
the equality $\calD(T)=\calD(T^+)$ holds if $V\in N^1_\kappa(R)$. In
this case we may not consider the bi-extensions of $T$ in the rigged
Pontryagin space and the corresponding l.s.c.d.s. can be written as
follows
\[
\theta =\begin{pmatrix}
T &K &\calJ\\
\Pk&{ } &E
\end{pmatrix}.
\]
\end{remark}

\begin{theorem} \label{t17}
Let $ \theta $ be a l.s.c.d.s. of the form {\rm(\ref{e3-2})} such
that $\calJ=I$, $\ol{\calD(\dA)}\ne\Pk$, and
$\calD(T)\ne\calD(T^+)$. Then operator-function $V_\theta (z)$ of
the form \eqref{e3-5}, \eqref{e3-6} has a holomorphic continuation
 $V(z)$ which belongs to the class $N^{01}_\kappa(R)$.

Conversely, let an operator-valued function $V(z)$  belongs to the
class $N^{01}_\kappa(R)$.  Then $V(z)$ admits a minimal
realization by a system $\theta$ of the form \eqref{e3-2} with a
non-densely defined operator $\dA$ and $\calD(T)\ne\calD(T^+)$.
\end{theorem}
The proof is immediate from Theorem~\ref{t10} and Theorem~\ref{t15}.

\section{Applications to the scalar case}
\label{s7}

In this section we consider scalar functions of the class
$N_\kappa$. We will establish the link between scalar ($E=\bbC$)
realizable functions of the class $N_\kappa(R)$ and a class of
realizable Nevanlinna functions \cite{BT3}.

The realization problems of the present type for Nevanlinna
operator-valued functions were studied in details in \cite{BT3}
and \cite{BT4} where similar subclass structure was developed. In
particular, it was shown that any realizable Nevanlinna
operator-function $V(z)$ admits an integral representation
\begin{equation}\label{e-62}
    V(z)=Q+F\cdot z+\int\limits_{-\infty }^{+\infty }\left( \frac 1{t-z}- \frac
{t}{1+t^2}\right)\,dG(t),
\end{equation}
in the Hilbert space $E$. In this representation $Q=Q^*$, $F=0$,
$G(t)$ is non-decreasing operator-function on $(-\infty,+\infty)$
for which
$$\int\limits_{-\infty}^{+\infty }\frac
{dG(t)}{1 + t^2} \in [E,E],$$ and
$$Qe = \int\limits_{-\infty} ^ {+\infty }
\frac {t} {1+t^2}\,dG(t) e$$ for all $e\in E_\infty$ where
\begin{equation}\label{e-63}
E_\infty=\left\{ e\in E:  \int_{-\infty }^{+\infty }\left(
dG(t)e,e\right) _E< \infty \right\}.
\end{equation}
 The class of all realizable
Nevanlinna operator-functions is called $N(R)$ (see \cite{BT3}).
The three subclasses of the class $N(R)$ were introduced in
\cite{BT4} and are called $N^0(R)$, $N^1(R)$, and $N^{01}(R)$.
Each subclass is described in terms of the subspace $E_\infty$ in
\eqref{e-63} determined by the variation of measure $G(t)$ in the
representation \eqref{e-62}. In particular, $E_\infty=\{0\}$ for
the class $N^0(R)$, $E_\infty=E$ for $N^1(R)$, and
$\{0\}\varsubsetneq E_\infty\varsubsetneq E$ for $N^{01}(R)$.

Now let us recall the following factorization result from
\cite{DLLSh} (see also~\cite{DHS}). Every scalar function $V(z)\in
N_\kappa$ admits a unique factorization
\begin{equation}\label{fact2}
    V(z)=\frac{p(z)p^\sharp(z)}{q(z)q^\sharp(z)}V_0(z),
\end{equation}
where $V_0$ belongs to the class $N_0$, $p(z)$ and $ q(z)$ are
relatively prime monic polynomials such that $\max(\mbox{deg }p,
\mbox{deg }q)=\kappa$, $p^\sharp(z)=\overline{p(\bar z)}$. Recall
also that the point $\infty$ is called a generalized pole of
nonpositive type of $Q$  if
\begin{equation}
\label{infgpol} -\infty\leq\lim_{z\widehat{\rightarrow }\infty }
\frac{Q(z)}{z} < 0,
\end{equation}
and the point $\infty$ is called a generalized zero of nonpositive
type  of $Q$ if
\begin{equation}
\label{infgzer}  0\leq \lim_{z\widehat{\rightarrow }\infty } {z}Q(z)
< \infty.
\end{equation}
In terms of the factorization~\eqref{fact2} one can consider the
following three possibilities:
\begin{enumerate}
\item
$\infty$ is a generalized pole of non-positive type  of the function
$V(z)$ if and only if $\mbox{deg }p > \mbox{deg }q$;
\item
$\infty$ is a generalized zero of non-positive type  of the function
$V(z)$ if and only if $\mbox{deg }p < \mbox{deg }q$;
\item
$\infty$ is neither the generalized pole of non-positive type nor
generalized zero of non-positive type of the function $V(z)$ if and
only if $\mbox{deg }p = \mbox{deg }q$.
\end{enumerate}
In the first case $V$ does not belong to the class $N_\kappa(R)$, in
the second case $V$ is definitely in the class $N_\kappa(R)$, and in
the third case the inclusion  $V\in N_\kappa(R)$ can be
characterized in terms of the function $V_0$.
\begin{theorem} \label{t19}
Let $V(z)\in N_\kappa$ with $E=\bbC$ and let $\infty$ be neither the
generalized pole nor generalized zero of non-positive type of the
function $V(z)$. Then $V(z)$ belongs to the class $N_\kappa(R)$ if
and only if the function $V_0$ in the factorization \eqref{fact2}
belongs to the class $N(R)$.
\end{theorem}
\begin{proof} Suppose $V_0(z)\in N(R)$ and let $V$ admits
the factorization \eqref{fact2} with $p$ and $q$ such that
$\mbox{deg }p = \mbox{deg }q=\kappa$. Then $V$
 belongs
to the class $N_\kappa$ (see~\cite{DLLSh}). Also, since $V_0(z)\in
N(R)$ then $F=0$ in \eqref{e-62}. It follows from the representation
\eqref{e-62} (see \cite{KK74}, \cite{GzT}) that in this case
\begin{equation}\label{e-66}
    \lim_{y\uparrow\infty}\frac{V_0(iy)}{y}=F=0.
\end{equation}
Combining \eqref{fact2} and \eqref{e-66} we get \eqref{e6-57}. In
order to prove the third item in the definition of the class
$N_\kappa(R)$ we first notice that it is  equivalent to the function
$V(z)$ not being an identical constant. Let us assume the contrary,
i.e.,
$$    V(z)=\frac{p(z)p^\sharp(z)}{q(z)q^\sharp(z)}V_0(z)\equiv
k,\quad k\in\bbC.$$ This immediately contradicts that $V_0(z)$ is
holomorphic in the upper half-plane. Therefore, the condition (3) of
the definition of the class $N_\kappa(R)$ is satisfied.

It was shown in \cite{KK74} that if $V_0(z)\in N(R)$ then
\begin{equation}\label{e-67}
    \lim_{y\uparrow\infty}yV_0(iy)=\int^\infty_{-\infty}dG(t),
\end{equation}
where $G(t)$ is the function from the representation \eqref{e-62} of
$V_0(z)$. Then we see that the subspace $\calB$ defined by
\eqref{e6-59} for the class $N_\kappa(R)$ coincides with the
definition of the subspace $E_\infty$ in \eqref{e-63} written for
the function $V_0(z)$. We notice that since $E=\bbC$ then either
$E_\infty=\{0\}$ or $E_\infty=E=\bbC$. Finally, if $E_\infty=\bbC$
$$\lim_{y\uparrow\infty}V(iy)=\lim_{y\uparrow\infty}V_0(iy)=0,$$
(see \cite{KK74}) and hence $V(z)\in N_\kappa(R)$.

Conversely, let $V(z)\in N_\kappa(R)$. Then \eqref{e-66} will
provide us with $F=0$ in the integral representation \eqref{e-62} of
$V_0(z)$. Furthermore, since $V(z)\in N_\kappa(R)$ then
\begin{equation}\label{e-68}
    \lim_{y\uparrow\infty}yV(iy)=\lim_{y\uparrow\infty}yV_0(iy),
\end{equation}
and is either finite or infinite. If the limit is infinite then
$E_\infty=\{0\}$ for $V_0(z)$ and $V_0(z)\in N(R)$. If the limit is
finite then $E_\infty=E=\bbC$, \eqref{e-67} holds for $V_0(z)$, and
(see \cite{KK74})
$$Q = \int\limits_{-\infty} ^ {+\infty }
\frac {t} {1+t^2}\,dG(t).$$ Therefore, $V_0(z)\in N(R)$.
\end{proof}
\begin{corollary} \label{t20}
A function $V(z)$  belongs to the class $N^{0}_\kappa(R)$ with
$E=\bbC$ if and only if the function $V_0$ in the equation
\eqref{fact2} belongs to the class $N^0(R)$.
\end{corollary}
\begin{corollary} \label{t21}
A function $V(z)$  belongs to the class $N^{1}_\kappa(R)$ with
$E=\bbC$ if and only if the function $V_0$ in the equation
\eqref{fact2} belongs to the class $N^1(R)$.
\end{corollary}
\begin{proof} The proofs of both corollaries immediately follow
from \eqref{e-67} and \eqref{e-68}.
\end{proof}

\section{Examples}\label{s9}
 We conclude the paper with simple illustrations.

\begin{example}
 Let us define $\Pi_1$ as a set of all $L^2({[0,2\pi]},dx)$
functions with the scalar product
\begin{equation*}
    [f,g]=\int\limits^{2\pi}_0f(x)\overline{g(x)}\,dx-\frac{1}{\pi}\int\limits^{2\pi}_0f(x)\,dx
    \int\limits^{2\pi}_0\overline{g(x)}\,dx
\end{equation*}
Let also $\dA$ be a $\pi$-symmetric operator defined by
\[
    \dA\,f=\frac{1}{i}\frac{df}{dx}
\]
with
\[
    \DdA=\Big\{f\in\Pi_1\mid f,f'\in
    AC_{\loc}([0,2\pi]),\,f(0)=f(2\pi)=0\Big\}.
\]
Let $T$ be an operator in $\Pi_1$ defined by
\[
    T\,f=\frac{1}{i}\frac{df}{dx}-\frac{1}{\pi i}\,f(2\pi),
\]
where
\[
    \calD(T)=\Big\{f\in\Pi_1\mid f,f'\in
    AC_{\loc}([0,2\pi]),\,f(0)=0\Big\}.
\]
One can check that the operator $T\supset\dA$, $T^+\supset\dA$,
and $\dA$ is a maximal $\pi$-symmetric part of $T$ and $T^+$,
i.e., $T\in\Omega_A$. The following two formulas define a
$(*)$-extension of $T$.
\[
\begin{split}
  &\bbA f= \frac{1}{i}\frac{df}{dx}-\frac{1}{\pi
  i}\big(f(2\pi)-f(0)\big)
  -if(0)\left[\delta(x-2\pi)+\delta(x)+\frac{2}{\pi}\right], \\
  &\bbA^\times f= \frac{1}{i}\frac{df}{dx}-\frac{1}{\pi i}\big(f(2\pi)-f(0)\big)
  +if(2\pi)\left[\delta(x-2\pi)+\delta(x)-\frac{2}{\pi}\right].
\end{split}
\]
By straightforward calculations we get
$$
\frac{\bbA-\bbA^\times}{i}f=\big(f(0)+f(2\pi)\big)\left[\frac{2}{\pi}-
\delta(x-2\pi)-\delta(x)\right].
$$
Now we can include $\bbA$ into a l.s.c.d.s.
\[
\theta =\left(\begin{array}{ccc}
\bbA &K &1\\
\mh^+\subset \Pi_1\subset \mh_- &{ } &\bbC
\end{array}\right),
\]
where
\[
    Ke=\frac{1}{\sqrt2}\left[\frac{2}{\pi}-
\delta(x-2\pi)-\delta(x)\right]e,\quad e\in\bbC.
\]
Then we can derive
\[
    W_\theta(z)=1+2i[(\bbA-zI)^{-1}Ke,Ke]=\frac{(z\pi i-1)e^{2\pi zi}+1}{e^{2\pi zi}-z\pi
 i-1}.
\]
Consequently, the function
\[
    V_\theta(z)=i [W_\theta (z) + I]^{-1} [W_\theta (z) - I] \calJ=\frac{2+\pi i z-(2-\pi i z)
    e^{2\pi i z}}{\pi z(e^{2\pi zi}-1)},
\]
belongs to the class $N_1^0$.
\end{example}

\begin{example}
Now we consider a construction which leads to examples of
functions of the class $N_\kappa^1(R).$ Let $\calH$ and $\calN$ be
two Hilbert spaces. Suppose that $A_0$ is possibly unbounded
operator in $\calH$ with nonempty resolvent set $\rho(A_0)$. Let
the operator $T$ in the Hilbert space $\calH=\calH\oplus\calN$ is
given by the block-operator matrix
\[
T=\begin{pmatrix}A_0&C\cr B&D\end{pmatrix}
\]
with bounded entries $B,C$ and $D$. Recall that the operator
valued function
\[
S(z)=D-B(A_0-zI)^{-1}C,\; z\in\rho(A_0)
\]
is called the transfer function of the system determined by the
matrix $T$. It is well known that the number $z\in\rho(A_0)$
belongs to $\rho(T)$ if and only if
\begin{equation}
\label{X}
 X(z):=S(z)-zI=D-B(A_0-zI)^{-1}C-zI
\end{equation}
has bounded inverse defined on $\calN.$  If this is a case then
the resolvent $(T-zI)^{-1}$ is given by the Schur-Frobenius
formula
\begin{equation}
\label{RT}
\begin{pmatrix}
(A_{0}-z I)^{-1}\left(I+CX^{-1}(z)B(A_0-z I)^{-1}\right)
&-(A_{0}-z I)^{-1}CX^{-1}(z) \cr - X^{-1}(z)B(A_{0}-z I)^{-1}
&X^{-1}(z)\end{pmatrix}
\end{equation}
Note that $X(z)$ is called the Schur complement of the
block-matrix $T-zI$. Let $\dim\calN=\kappa<\infty$. Suppose that
$A_0$ is a selfadjoint operator in $\calH$. In this case if $|\Im
z|$ is sufficiently large then the norm of the operator
$z^{-1}\left(D-B(A_0-zI)^{-1}C\right)$ is less then one, therefore
$X^{-1}(z)$ exists and $X^{-1}(z)$ is a meromorphic function in
$\bbC_+\cup\bbC_-$.  Equip the Hilbert space $\calH$ by the
indefinite inner product
\[
[h,g]:=(P_{\calH_0}h,P_{\calH_0}g)- (P_{\calN}h P_{\calN}g),\;
h,g\in\calH,
\]
where $P_{\calH_0}$ and $P_{\calN}$ are the orthogonal projections
in $\calH$ onto $\calH_0$ and $\calN$, respectively. Then $\calH$
becomes a Pontryagin space $\Pi_\kappa.$ Let $A$ be a linear
operator in the Hilbert space $\calH$ defined as follows
\[
\calD(A)=\calD(A_0),\;Ah:=A_0h+Bh,\; h\in\calD(A).
\]
The operator $A$ is a non-densely defined and closed $\pi$-Hermitian
operator in $\Pi_\kappa$. Since the operator $P_{\calH_0}A=A_0$ is
selfadjoint, the operator
 $A$ is $J$ regular $O$-operator, where $J=P_{\calH_0}-P_\calN$. Evidently, the operator
\begin{equation}
\label{OpT}
 T=\begin{pmatrix}A_0&-B^*\cr B&D\end{pmatrix}
\end{equation}
meets the conditions
\[
\calD(T)=\calD(T^+)=\calD(A_0)\oplus\calN,\;T\supset A,\;
T^+\supset A
\]
In particular, $T$ is $\pi$-selfadjoint if and only if $D$ is a
selfadjoint in the Hilbert space $\calN$.

Let $D=\Re(D)+i\Im(D)$. Then
\[
\Re(T)=\frac{1}{2}(T+T^+)=\begin{pmatrix}A_0&-B^*\cr
B&\Re(D)\end{pmatrix},\;
\Im(T)=\frac{1}{2i}(T-T^+)=\begin{pmatrix}0&0\cr 0
&\Im(D)\end{pmatrix}
\]
 Assume also that $D_I=K\calJ
K^+$, where $K$ acts from the Hilbert space $\calN$ into the
negative subspace $\calN$ of the Pontryagin space $\Pi_\kappa$  and
$K$ is invertible. Note that $K^+=-K^*$ where $K^*$ is the Hilbert
space adjoint to $K:\calN\to\calN$. Then
\[
\theta =\begin{pmatrix}
T &K &\calJ\\
\Pk&{ } &\calN
\end{pmatrix}
\]
is the l.s.c.d.s. The transfer function
 $W(z)$ of $\theta$ is given by
\[
W(z)=I-2iK^+(T-zI)^{-1}K\calJ=I-2iK^+P_\calN(T-zI)^{-1}K\calJ,\;z\in\rho(T)
\]
and its fractional-linear transformation is
\[
V(z)=i[W(z)+I]^{-1}[W(z)-I]\calJ=K^+P_\calN (\Re(T)-zI)^{-1}K,\;
z\in\rho(T)\cap\rho(\Re(T)).
\]
 Let
\[
X_T(z)=D+B(A_0-zI)^{-1}B^*-zI,\;
X_{T_R}(z)=\Re(D)+B(A_0-zI)^{-1}B^*-zI.
\]
Then from \eqref{RT} it follows that
\[
W(z)=I-2iK^+X^{-1}_T(z)K\calJ,\; V(z)=K^+X^{-1}_{T_R}(z)K.
\]
Let us take $\Im(D)=-I$, $K=I$ then $K^+=-I$ and $-X^{-1}_{T_R}(z)$
belongs to the class $N_\kappa^1(R)$ in the Hilbert space $\calN.$
Thus we obtain the following theorem.

\begin{theorem}
\label{NR} Let $\calH_0$ and $\calN$ be Hilbert spaces,
$\dim\calN=\kappa<\infty$.  Let $A_0$ be a selfadjoint operator in
$\calH_0$, $B$ is a bounded operator from $\calH_0$ into $\calN$,
and let $D$ be a selfadjoint operator in $\calN$. Then the
operator valued function
\[
V(z)=-\left(D+B(A_0-zI)^{-1}B^*-zI\right)^{-1}, \; z\in
\bbC_+\cup\bbC_-
\]
belongs to the class $N_\kappa^1(R)$ in the Hilbert space $\calN$.
\end{theorem}
Using Theorem \ref{NR} let us give some concrete example of scalar
function from the class $N_1^1(R)$.

 Let $\calH$ be a weighted Hilbert space
$\calL_2\left([-1,1],\rho(t)\right)$ with the weight
\[
\rho(t)=\frac{2}{\pi}\, \sqrt{1-t^2}.
\]
Let the operator $A_0$ in $\calL_2\left([-1,1],\rho(t)\right)$ be
defined as follows:
\[
(A_0f)(t)=tf(t),\; f(t)\in \calL_2\left([-1,1],\rho(t)\right).
\]
Then $A_0$ is a selfadjoint contraction. Let $e_0(t)=1$, $t\in
[-1,1]$ The function $e_0(t)$ belongs to
$\calL_2\left([-1,1],\rho(t)\right)$ and $||e_0||=1.$  Let
$\calN=\bbC$. Define the operator
$B:\calL_2\left([-1,1],\rho(t)\right)\to\bbC $ as follows
\[
Bf(t)=\frac{2\gamma}{\pi}\,\int\limits_{-1}^1
f(t){\sqrt{1-t^2}\,dt},\;f(t)\in \calL_2\left([-1,1],\rho(t)\right),
\]
where $\gamma\ne 0$. Then
\[
B^* c=\overline\gamma\,c\,e_0(t),\;c\in\bbC.
\]
Let $D$ be the operator of multiplication on a real number $d$ in
the space $\bbC$.  It is known  \cite{Ber} that
\[
\frac{2}{\pi}\int\limits_{-1}^1\frac{\sqrt{1-t^2}}{t-z}\,dt=2(\sqrt{z^2-1}-z),\;
z\notin [-1,1],
\]
where the branch of the function $\sqrt{z^2-1}$ is taken such that
$\Im\,\sqrt{z^2-1}>0 $ for $\Im\,z>0$. It follows that the function
$V(z)=-\left(D+B(A_0-zI)^{-1}B^*-zI\right)^{-1}$ takes the form
\[
V(z)=\frac{1}{z-2\,|\gamma|^2(\sqrt{z^2-1}-z)-d}=\frac{1}{(1+2|\gamma|^2)z-2\,|\gamma|^2\sqrt{z^2-1}-d}
\]
According to Theorem \ref{NR} the function $V(z),\; z\in \bbC_+$
belongs to the class $N_1^1(R)$. If
$|\gamma|^2>\max\{0,\,(d^2-1)/4\}$ the function $V(z)$ has a simple
pole
\[
z_0=\frac{d(1+2|\gamma|^2)+2|\gamma|^2i\sqrt{4\,|\gamma|^2+1-d^2}}{1+2|\gamma|^2}
\]
in $\bbC_+.$

\end{example}


\end{document}